# SIMULATED ANNEALING: IN MATHEMATICAL GLOBAL OPTIMISATION COMPUTATIONS, HYBRID WITH LOCAL OR GLOBAL SEARCH, AND PRACTICAL APPLICATIONS IN CRYSTALLOGRAPHY & MOLECULAR MODELLING OF PRION AMYLOID FIBRILS

*Jiapu Zhang*[*]

School of Science, Informatics Technology, and Engineering, &
Centre of Informatics and Applied Optimisation,
The University of Ballarat, MT Helen Campus, Victoria 3353, Australia

## Abstract

Simulated annealing (SA) was inspired from annealing in metallurgy, a technique involving heating and controlled cooling of a material to increase the size of its crystals and reduce their defects, both are attributes of the material that depend on its thermodynamic free energy. In this Paper, firstly we will study SA in details on its initial feasible solution choosing, initial temperature selecting, neighborhood solution searching, efficient way of calculating for the difference of objective function values of two neighborhood solutions, acceptance function (Metropolis function), temperature cooling, and the criteria of inner and outer loops' stopping, etc. Then, hybrid pure SA with local (or global) search optimization methods allows us to be able to design several effective and efficient global search optimization methods. In order to keep the original sense of SA, we clarify our understandings of SA in crystallography and molecular modeling field through the studies of prion amyloid fibrils.

**Keywords**: Simulated Annealing; Computational Algorithm; Global Optimization Search; Hybrid with Local/Global Search; Applications to Crystallography.

## 1. Introduction

Simulated Annealing (SA), as well as Tabu Search and Genetic Algorithm, is one of the successful heuristic computational methods. It simulates the annealing process with Monte Carlo property. The works of Metropolis, Kirkpatrick, Johnson, Aarts, et al. are well-known. In discrete optimization, simulated annealing method has found a lot of applications. The book for example [1] is good collections of its applications to discrete optimization problems. For continuous optimization problems, there are a lot of references on it. However, we still rarely see a very successful simulated annealing method for large scale continuous optimization problems in very high dimensions, especially in the constrained case.

The SA method appeared as early as in 1953 [2] as a Monte Carlo method and was firstly investigated and used in 1983 by Kirkpatrick et al. [3]. SA is a stochastic method. It differs from the traditional descent methods (see, for example, [4] and references therein) in that a local search algorithm for a neighborhood solution search, whether it randomly descents or steeply descents, allows downhill moves only, while in an attempt to escape local optima SA algorithm allows occasional uphill moves as well. SA techniques are based upon the physical analogy of cooling crystal structures (including the case of so called quenching) that spontaneously arrives at a stable configuration, characterized by globally or locally-minimal potential energy. Starting with an initial solution x, and an initial "temperature" $T_0$, which is a parameter, we obtain a neighboring solution x′ and compare its cost with that of x. If the cost of x′ is smaller than that of x, i.e. f(x′)< f(x), we accept the new solution x′. The same thing would happen if we are applying the local search method by random descent method [4]. On the other hand, if f(x′) is greater than f(x), (in which case local search algorithms (see, for example, [4]) will not accept x′), the SA algorithm may accept x′, but with a probability $e^{-\Delta x'x/T_0}$ where $\Delta x'x$ is the difference in the costs of x′ and x, i.e., $\Delta x'x = f(x′)−f(x)$. This process is carried out for a certain number of times, which we call iterations, for each temperature. Then we reduce the temperature according to a particular schedule, and repeat. The convergence of SA algorithms are studied, for example, in [5-6].

---

[*] E-mail address: j.zhang@ballarat.edu.au, jiapu_zhang@hotmail.com, Phone: 61-4 2348 7360, 61-3 5327 6335

An essential element of the SA algorithm is the probability $e^{-\Delta x'x/T}$ of an uphill move of size $\Delta_{x'x}$ being accepted when the current temperature is T. This is dependent on both $\Delta x'x$ and T. For a fixed temperature T, smaller uphill moves $\Delta_{x'x}$ have a higher probability of being accepted. On the other hand, for a particular uphill move $\Delta_{x'x}$, a higher temperature will result in a larger probability for that uphill move being accepted. As stated in [7], "at a high temperature any uphill move might be indiscriminately accepted with large probability so that the objective function and the tumbles around the space are not very important; as T is lowered the objective function becomes more and more significant; until as T goes to zero the search becomes trapped in the lowest minima that it has reached."

The SA algorithm for solving a practical problem is typically implemented in two nested loops: the outer loop and the inner loop. The outer loop controls temperatures, while the inner loop iterates a fixed number of times for the given temperature. The inner loop is for the problem specific decisions. The decisions of the outer loop involve the setting of initial temperature ($T_0$), the cooling schedule, the temperature length which is the number of outer loop iterations performed at each temperature, and the stopping condition of the outer loop. The inner loop of SA typically considers the following aspects: feasible solution space, initial feasible solution, neighborhood move, objective function values (and efficient calculation of their difference), and the decision which decides whether the move is found acceptable or probably acceptable according the so-called Metropolis criterion.

In this paper, Section 2 will study SA in details on its initial feasible solution choosing, initial temperature selecting, neighborhood solution searching, efficient way of calculating for the difference of objective function values of two neighborhood solutions, acceptance function (Metropolis function), temperature cooling, and the criteria of inner and outer loops' stopping, etc. In Section 3, hybrid pure SA with local (or global) search optimization methods allows us to be able to design several effective and efficient global search optimization methods. In order to keep the original sense of SA, we clarify our understandings of SA in crystallography and molecular modeling field through the studies of prion amyloid fibrils in Section 4. Section 5 will give some concluding remarks on SA.

## 2. Implementing the SA algorithm

### 2.1 Overview

In this section we consider problem

$$\text{Minimize } f(x) \quad \text{subject to } x \in X,$$

where X is a subset of $R^n$ is a compact set and f is continuous, being solved by SA algorithm.

The word *renew* denotes the counts of the solution being accepted in the inner loop. The pseudo-code (referred, for example, to [8]) of the SA algorithm is listed as follows:

Algorithm 1. *Simulated annealing algorithm*.

**Initialization**:
    Define the objective function f
    and its feasible solution space.
    Call **Initial feasible solution producing procedure**
    to produce an initial feasible solution x.
    Call **Procedure of selecting initial temperature**
    to produce the initial temperature $T_0$.
    Calculate the size of neighborhood N_size.
    Calculate f(x), and set x_best = x and f_best = f(x).
    Set best_count = frozen_count = 0, and value of δ.

**Cooling (outer loop procedure)**:
    **Repeat** (outer loop)
        Call **Inner loop procedure**.
        Call **cooling schedule** T = α(T)
        to decrease to a new temperature.
        **If** best_count > 0 **then** set frozen_count = 0
        **If** renew/iteration count < 1/N_size **then**
            set frozen_count = frozen_count + 1
    **Until outer loop stopping criterion** is met

**Inner loop Procedure**:
Set iteration_count = 0.
**Repeat** (inner loop)
  Call **Neighborhood solution search procedure**
  to generate a feasible neighborhood solution x′.
  Calculate f(x′).
  Call **Efficient procedure**
  of calculating the cost difference $\Delta_{x'x} = f(x') - f(x)$.
  **If** $\Delta_{x'x} < -\delta$ **then**
    Set x = x′, renew = renew + 1.
    Set f(x) = f(x′).
    **If** f(x) < f_best **then**
      x_best = x
      f_best = f(x)
      best_count = best_count+1
      Record results on "Best So Far"
    **Endif**
  **else**
    **If** random[0, 1] < exp($-\Delta_{x'x}$/T) **then**
      x = x′
      f(x) = f(x′)
      renew = renew + 1
    **Endif**
    Set iteration_count = iteration_count + 1.
  **Endif**
**Until stopping criteria of inner loop** is met

In implementing the SA algorithm described above, *initial feasible solution producing procedure, the procedure of selecting initial temperature, neighbourhood solution search procedure, efficient way to calculate the difference of objective function values of two neighbourhood solutions, acceptance function (here it is the Metropolis function), cooling scheduling of the temperature, and stopping criterions of inner* and *outer loops* are its important components. Different definitions of those are discussed in the literature on SA methods, which will be discussed in the following subsections.

SA algorithm corresponds to a Markov chain. For each temperature T fixed, if the variation of Markov chain arrives at a stable state and then T goes down, we call the SA algorithm *homogeneous SA algorithm*; if not, it is a *non-homogeneous SA algorithm*.

### 2.2 Initial feasible solution producing procedure

For a convex function the initial feasible solution can be chosen anywhere, from which the global minimum is reached by moving towards to the lower values, in the feasible region; however, for a non-convex function, it depends on the initial solution very much either to find a local minimum of it or to find its global minimum [9]. For real projects, usually there are many requirements, i.e. constraints, for reaching its aims. A good feasible initial solution producing procedure is clearly needed. Numerical experiments show that, without the sensitive procedure of choosing the initial simplex (see § 3.1.1 of [10]), the Simplex Simulated Annealing (SSA) method of paper [10] is very difficult to make it work. However, for many problems solved by SA method, there, often, is a simple way of producing initial solution: randomly taking a feasible solution from the feasible region as the initial solution. Constraint programming is a new high-level language paradigm for satisfaction and optimization problems. To produce a feasible solution by constraint programming strategy as the initial solution is also a very popular way; for example, see [11]. Using a local/global search optimization method to quickly get a solution as the initial solution of SA is also a good way.

### 2.3 Initial temperature selecting procedure

Numerical experiments tell us that proper initial temperature $T_0$ can make the SA method quickly get the optimal value of the objective function. If at initial temperature we accept almost all the solution (i.e. acceptance rate $\chi_0 \approx 1$), then in theory, by Metropolis criterion exp $\Delta_{x'x}/T_0 \approx 1$, where $T_0$ should be "sufficiently" large. Johnson, Kirkpatrick, Aarts et al. [12-15, 3] present several initial temperature selecting procedures. The idea of Kirkpatrick is: to choose a large $T_0$, give $\chi_0$ in advance (for example $\chi_0 = 0.8$), generate many solutions, if the acceptance rate $\chi$ is less than $\chi_0$ then increase $T_0$, repeat until $\chi > \chi_0$ to get a $T_0$. Johnson's formula is
$$T_0 = \Delta^- f^+ / \ln(\chi_0^{-1}),$$
where $\Delta^- f^+$ is the average increase of objective function values of randomly generated solutions. The one of [12] is frequently used; for example, in [10, 8].
Aarts' formula is
$$T_0 = \Delta^- f^+ / \ln ( m_2 / ( m_2 \chi - m_1 (1-\chi) ) ),$$
where $m_1$ is the number of solutions making the objective function value decrease, $m_2$ is the number of solutions making the objective function value increase, and $\chi$, for example, may be set as $\chi_0$. However, those procedures are not definitely working well for all problems. Fixed temperature schedule is studied and applied in [16-18]. In homogeneous SA, $T_0$ chosen should be properly large enough to sufficiently accept all the candidate feasible solutions possibly produced. We also often use the following ways:

(1). Uniformly sample some solutions, calculate their objective function values, and take the variance of those objective function values as $T_0$.

(2). Randomly generate some solutions, determine $|\Delta_{max}|$ which is the maximal difference of each pair of solutions, and calculate $T_0 = -\Delta_{max} / \ln p_r$, where pr is the initial acceptance probability and in theory it should be close to 1. For non-homogeneous SA, in theory we have formulas [19] for calculating $T_0$.

**2.4 Neighborhood solution searching procedure**

This is one key element in implementing SA. For discrete optimization problem, for instance, in the book [1], there are several successful neighborhood solution searching procedures. All those procedures should be at least based on two basic ideas: (a) neighbor means "nearby", (b) SA method is a stochastic method so that the neighborhood solution should be randomly taken. We may take those ideas for developing neighborhood solution searching procedure for continuous optimization problems.

First we review some neighborhood solution searching procedures of continuous optimization problems. Miki et al. (1999) presented a formula $x^{k+1}_i = x^k_i + r*m$, where r is a random number with uniform distribution in $[-1,1]$, m is the neighborhood range which makes the rate between accepted and rejected moves approximate 0.5 [20]. In [21], first generate a random direction vector $\theta_k \in R^n$, with $||\theta_k|| = 1$, then find a fixed step size $\Delta r$, thus get a neighborhood solution $x_{k+1}$ of $x_k$: $x_{k+1}=x_k+\Delta_r\theta_k$. The choice of $\Delta r$ is thoroughly discussed in [22-23]. In [23], the direction vector $\theta_k$ is defined in a new way. It is suggested to take into account the point $x_h$, h<k, generated by the algorithm and different from $x_k$, if $f(x_h)<f(x_k)$ then $\theta_k=x_h-x_k$, otherwise $\theta_k=x_k-x_h$. Contrary to the results in [21-22], Corana et al. (1987), Siarry et al. (1997) and Vanderbilt et al. (1984) search through the space of feasible region in an anisotropic way [24-26]. In [24], at each iteration k a single variable of $x_k$ is modified, and iterations are subdivided into cycles of n iterations during which each variable is modified; i.e. $x_{k+1}=x_k + r\ v_{i+1}\ e_{i+1}$, where r is a uniform random number in [-1,1], $i\in\{0, …, n-1\}$ is such that k+1 = h*n+i for some nonnegative integer h, and $v_{i+1}$ (that is anisotropic) is the maximum allowed step along the direction $e_{i+1}$ of the (i+1)-st axis. Instead of varying a single variable in $x_k$ in each iteration, Siarry et al. (1997) varied p variables [25]. Another concept to simulated annealing method is adaptive (see [27-32]). This means SA method should possess the ability of adapting itself to the problem it solves, the objective function f and the temperature, etc. whether globally or locally. The code of Ingber's ASA (Adaptive Simulated Annealing) algorithm [28-31] can be retrieved from the web site www.ingber.com, and many techniques such as 'fast annealing', 're-annealing', 'quenching', 'multistart strategy', and 'priori information' are used. Romeijn et al. (1994, 1999) proposed a two-phase generator: "*first generating a random direction $\theta_k$, with $||\theta_k|| = 1$, and then generating a random point $\lambda_k$ in the set $\Lambda_k= \Lambda_k(\theta_k)=\{\lambda: x_k+\lambda\theta_k \in X\}$, thus $x_{k+1}=x_k+\lambda_k\theta_k$*", and if $x_{k+1}$ is not in X or if there is a jamming problem, i.e. $\Lambda_k$ is very small, then use the '*reflections*' technique [33-34]. Employing computer science theory is also useful for the neighborhood solution searching procedure; for example Bilbro et al. (1991) gave a tree annealing approach [35]: *divide the feasible regions in the form of a tree, and $x_{k+1}$ is sampled from a distribution which selects a random path from the root of the tree to a leaf of the tree in a way that the sub-regions with a high percentage of accepted points are favored*. Employing a local search method into simulated annealing method is also very popular. Desai et al. (1996) proposed a technique [36]: *randomly perturb the current solution $x_k$ to get a new point $x_{k+1}$, and start a local search from $x_{k+1}$ to get a new local minimum $x_{k+1}$*, which attempts to combine the robustness of annealing in rugged terrain with the efficiency of local optimization methods in simple search spaces. The parallel version of [36] may be seen in [37]. Lucidi et al. (1989) presented a random tunneling technique by means of acceptance-rejection sampling [38]. Over the unit hypercube feasible region, Fox (1995) gave a special neighborhood solution searching procedure [39]. First the objective function f is evaluated at the points of a net. Then, the unit hypercube is subdivided into many boxes of a set C, which are with different widths along different axes. Over C, generate a probability mass function p by intersecting the net with each box to generate many points and find the minimum value of those points. Then sample a box B from C according to the probability mass function p, sample a uniform point y from B and apply some local search steps starting from y. Repeat the sampling of B and y a finite number of times, and get the set $F(x_k, k)$, which is a finite set of candidate neighborhood points of $x_k$ at iteration k. And then the acceptance probabilities are applied to define the distribution of the next iteration $L \geq k$, and according to the acceptance probabilities the algorithm randomly selects a point in $F(x_k, k)$ and moves to it. For getting neighborhood solution, there is an idea: simultaneously perturbing all of the variables of the problem in a "proper" random fashion. In [40] the neighborhood solution producing procedure is given by the way: *randomly uniformly re-generate one element of $x_k$, or $m \in random\{1, ...,n\}$ elements of $x_k$, or the whole vector of $x_k$, as the new solution $x_{k+1}$*. In [10], instead of the reflect-expand-contract-shrink Nelder-Mead method (seen, for example, in [41]), the reflects-shrink Simplex Direct Search (SDS) method is given. The SDS method uses various numbers of reflections in a flexible manner and follows a shrinking if after reflecting all the n worst vertices of the initial simplex better movement is still failed to be gotten.

In Fast Simulated Annealing method of [42], Cauchy distribution is used to produce new solution. Greenes et al. (1986) used the probability of fitness function, which is based on objective function, to produce new solution [43].

From the ideas of all those reviewed above, we present two versions of the neighborhood solution search procedure for SA algorithm. In the SDS algorithm, for the objective function value of each vertex we add a random fluctuation: $f(x_i)+k_B T\log(z)$, where $k_B$ is the Boltzmann constant in appropriate units and z is a random number of the interval (0,1). We might carry out a multiple shrinking, in which the highest vertex is simultaneously moved along all the coordinates towards the lowest-energy vertex. For getting a new neighborhood solution, in [40], the procedure below is used:
"*i = random{1, 2, 3}, which is a random integer taken from the set {1, 2, 3}. Depending on the outcome of i, within the feasible region, re-generate randomly one of the following: one element of x, or $m \in \{1, 2, \cdots, n\}$ elements of x, or the whole vector of x. This gives x'.*"

Noticing that the neighborhood solution search for simulated annealing method should be at least based on ideas (a) and (b), we may simply give a neighborhood solution search procedure for simulated annealing algorithm:

*Uniformly randomly keeping n−1 elements of x, and making the left one element of x uniformly randomly take a value such that the new solution x' is still feasible. This gives x'.*

When the feasible region of the optimization problem is the unit simplex S, the neighborhood solution search procedure should be modified:

*Uniformly randomly keeping n − 2 elements of x, and making one element from the two elements left to x uniformly randomly take a value from [0,1] such that the value of the sum of the n − 1 elements is not greater than 1. Another left element of x' is given the value 1-sum. This gives x'.*

This is an efficient procedure, which is better than the one in [40].

Another version of neighborhood solution procedure may be found in the following SA pseudo-code:

*New version of the simulated annealing algorithm.*

X_best = x & f_best = f; q := $q_0$
DO j := 1 to J
    T := $T_0(j)$
    Repeat
      DO k := 1 to L
        Repeat
            Randomly generate the search direction d ∈ (−1, 1)
            Let $x'_i = x_i + q \cdot d$ and $x'_l = x_l$ when $l \neq i$
        Until x' is feasible
        Calculate _x'x
        IF ($\Delta_{x'x}$ < 0) or (exp($-\Delta x'x/T$) > random[0, 1]) THEN
           Accept x'
           IF f < f_best THEN x_best = x & f_best = f
           Calculate next annealing temperature T
           Adjust the step length q := g(Acc) * q (where g(·) is a function given)
      END DO
    Until outer loop stopping criterion is satisfied
    IF f best < f THEN x = x best & f = f best
END DO

where g(·) is an adjustment function, for example g(x)=$(x-0.5)^3$+1. As a whole, the new solution generating procedure composes two parts: the way to generate candidate solution, and how to generate the probability distribution of the candidate solution. Hence, we may replace our uniform distribution by Cauchy distribution, Gauss normal distribution, or their combined distribution, and get some new results for comparisons.

## 2.5 Efficient calculation of cost difference

Since a very large number of iterations are performed, it is essential to calculate efficiently the cost differences between a solution and its neighborhood solution. Take a simple instance, if f(x)=Ax+b, it is clearly $\Delta_{x'x}$=A(x′−x) is much efficient than $\Delta_{x'x}$=(Ax′+b)−(Ax+b), especially when the computational effort is very much for the computer.

Note: There is a subtle difference in the meaning between objective function and fitness function. *Objective function* measures the variable's performance in the search region; whereas *fitness* function provides a measure of variable's relative performance by transforming the objective function f into F(x)=g(f(x)). For example, proportional fitness assignment ($F(x_i)=f(x_i)/\sum f(x_i)$), linear transformation (F(x)=af(x)+b) are some simple transformations.

## 2.6 Acceptance function

The (Markov chain state) acceptance function, generally, is given in a probability form that should meet the following criteria:
(1). At each fixed temperature, the acceptance function should maintain the average percentage of accepted moves at about 1/5 of the total number of moves, which cannot make the objective function value decrease;
(2). With the decrease of temperatures, the probability of accepting an increasing move decreases.
(3). When temperature becomes zero, only the solutions that make the objective function value decrease can be accepted.

In Subsection 2.1, the acceptance function for the SA algorithm is the so-called Metropolis function
$$A(x, x', T)=\min\{1, \exp\{-\Delta_{x'x}/T\}\}. \tag{2.1}$$
Note here we let the acceptance function depend on the difference of the function values of x and x′ instead of depending directly on x and x′. More generally, we may rewrite (2.1) as follows:
$$A(x, x', T)=\min\{1, \exp\{-\Delta_{x'x}/\gamma(T)\}\}, \tag{2.2}$$
if $\gamma$: $(0,+\infty)\to(0,+\infty)$ is a strictly increasing function under some balance conditions [44]. Barker's function
$$A(x, x', T)=1/(1+\exp\{-\Delta_{x'x}/T\}) \tag{2.3}$$
is another popular acceptance function. This function has a similar form for T to (2.2). Johnson et al. (1987) suggested
$$A(x, x', T)=\min\{1, 1-\Delta_{x'x}/T\}. \tag{2.4}$$
and made the speed of SA algorithm increase by 30% [45]. Sechen (1988) used table search to reduce the time wasted on calculating $\exp\{-\Delta_{x'x}/T\}$ [46].

## 2.7 Cooling scheduling of temperature

During the SA iterations, the temperature sequence $\{T_k\}$ is being produced. If $\lim_{k\to+\infty}T_k=0$, we say that $\{T_k\}$ is a *cooling schedule*. In this subsection, we review some successful cooling schedules.

Aarts and Laarhoven (1985) presented a cooling scheme [47]
$$T_{k+1}=T_k/(1+T_k\log(+\varepsilon)/(3\sigma_k)), \tag{2.5}$$
where $\sigma_k$ is the standard deviation of the observed value of the cost function, and $\varepsilon$ is 0.1 in [18].
$$T_{k+1}=T_k\exp(-T_k(f_{T_k} - f_{T_{k-1}})/\sigma^2_{T_k})$$
is another cooling schedule [48]. Reeves (1995) described a cooling schedule of Lundy and Mees [4], where the temperature is reduced according to
$$T_{k+1}=T_k/(1+\beta T_k), \tag{2.6}$$

or equivalently,

$$T_k = T_0/(1+k\beta T_0), \quad (2.7)$$

where $\beta$ is suitably small, and only one iteration is performed in each inner loop. For the convergence of non-homogeneous SA method, in 1984 Geman et al. (1984) gave the Boltzmann annealing or called classical SA [49], in which the temperature is calculated by

$$T_k = T_0/\ln(k+c), \quad k=1,\ldots,\infty, \quad (2.8)$$

where $c=1$. A little modification of $c$ is used in [50]:

$$T_k = T_0/\ln(k+c), \quad k=0,\ldots,\infty \quad (2.9)$$

with $c=e=2.7183$. For formulas (2.8) and (2.9) $c$ should not be less than 1. In 1987, Szu et al. (1987) proposed the Fast Annealing method [42]. The cooling schedule of this method is with a faster decrease:

$$T_k = T_0/(k+1), \quad k=1,\ldots,\infty \quad (2.10)$$

that decreases sharper than (2.8). However, we should match the rate of temperature decrease with the neighborhood solution generating procedure. Nahar et al. (1987) divided $[0, T_0]$ into K intervals and find $T_k$, $k = 1, \ldots, K$ [15]. The Very Fast Simulated Re-annealing method [28] was presented in 1989 by Ingber. Its cooling schedule is

$$T_k = T_0 \exp(-ck^{1/n}), \quad k=1,\ldots,\infty, \quad (2.11)$$

where $c$ is a scale factor. Ingber (1989) also used a slower schedule [28] of (2.10), which is

$$T_k = T_0/(k+1)^{1/n}, \quad k=1,\ldots,\infty. \quad (2.12)$$

Although many cooling schedules are mentioned above, the geometric cooling scheme [3] proposed by Kirkpatrick et al. (1983)

$$T_{k+1} = \alpha T_k, \quad k=0,\ldots,\infty, \quad (2.13)$$

where $\alpha \in (0, 1)$ is a constant, is still a widely used and a popular SA cooling schedule (refer to [12-15, 10, 4, 40]) because it compromises the quality and CPU time of optimization. Kirkpatrick et al. (1983) take $\alpha=0.95$; and Johnson et al. take $\alpha \in [0.5, 0.99]$. Our numerical experiments also shows that (0.8, 0.99) is a good interval chosen for $\alpha$. Given $T_0$, $T_f$ and the number of outer loop iterations, in [51] graphs of many kind of cooling schedules can be seen.

## 2.8 Stopping criterion of inner loops

The number of iterations in each inner loop is also called the temperature length. In many forms of simulated annealing method, a fixed number of iterations are performed for each temperature. Usually the fixed number is detected by a long sequence of iterations in which no new solutions have been accepted. This fixed number depends on the size of the neighborhood Nsize, which is defined to be the total number of possible distinct neighborhood solutions that can be obtained, and its mathematical form is Nfactor*Nsize, where Nfactor is some multiplying factor, for example Nfactor=10. In our pseudo-code of the simulated annealing method, we also introduce a symbol renew which records the number of times the solutions are accepted at a temperature. We may also terminate the inner loop if this number has exceeded Cut*Nfactor*Nsize where Cut is another multiplying factor. Section 4.2.3 of [9] describes this fixed number in view of stochastic process terminologies. In homogeneous SA, from the view of objective function values two stopping criteria may be presented:
(1). Checking stability of the expectation value of objective function values; and
(2). The change of objective function values is lower than some tolerance in a certain amount of iterations.

## 2.9 Stopping criterion of outer loops

The *choice of final temperature $T_f$* determines a stopping criterion of outer loops. At the end of each inner loop, if the best solution obtained in that inner loop has not been changed and at the same time we are not having many changes in the current solution, we reduce the temperature and start another inner loop. If a solution $x^*$ has been consecutively "frozen" at many current temperatures, then we stop and say $x^*$ is the best solution found by the simulated annealing method. Usually we set a proper small temperature as $T_f$ as the stopping condition of outer loops. Our numerical experiments show that we may get the small temperature T when "Floating point exception (core dumped)" is reminded by computer. In our pseudo-code, we also count the number *best_count* of times at which the best feasible solution is not replaced again, and calculate the proportion of solutions accepted, by *frozen_count*. We may halt the algorithm when *frozen_count* reaches a predetermined number. As a whole, based on $\lim_{k \to +\infty} T_k = 0$, we find a "frozen" temperature for the stopping criterion of outer loops. Nahar et al. (1987) used the number of temperatures, in other words, the number of Markov chains or iterations, as stopping criterion of outer loops [15]. Notice here the lengths of Markov chains, $L_k$, may be upper bounded by a constant $\bar{L}$. Johnson et al. (1987) used the acceptance rate to terminate the outer loops [45]: current acceptance rate $\chi_k > \chi_f$ given, where $\chi_f$ is the final acceptance rate.

From the point of view of the objective function values, we may also give the terminating criterion for the outer loops. If $|\Delta_{x'x}| \leq \varepsilon$ or $f(x) - f\_best \leq \varepsilon$, $|(f\_best - f\_opt)/f\_opt| \leq \varepsilon$ (where $|f\_best - f\_opt| < \varepsilon$ if $f\_opt = 0$), where $\varepsilon$ is a sufficiently small positive number and f_opt is the optimal value known, we stop the method. This simply means when the objective function values cannot be improved we may stop the algorithm.

In another form, we use the information not only on the temperature but also on the objective function values; then we can also give a stopping criterion for the simulated annealing method. Suppose $P_F$ is a proper number given, if $A(x, x', T) \leq P_F$, we stop the simulated annealing method. If in many successive Markov chains the solution has not changed, we can stop the method.

Setting an upper limit of executing time is also a way to stop the algorithm. The user may terminate the method manually according to a user-defined aim.

**2.10 Improvements on SA method**

Aarts et al. (1989) improved the simple cooling scheme of Johnson et al. talked in the above Subsection 2.7 (see [5]). They present a more meticulous cooling scheme in which $T_0$, $T_f$, $L_k$ and the formula of $T_k$ are well designed. Other improvements on SA are:
(1). Re-increase temperature. In the process of SA, in order to adjust some state, to reincrease its temperature is a good way.

*Heating-annealing procedure.*
$T_0 = 0$
Repeat
    DO k := 1 to L
        Generate new solution x′
        Calculate $\Delta_{x'x}$
        IF $\Delta_{x'x} > 0$ THEN accept x′ and *HEAT*: T := heat(T);
        IF ($\Delta_{x'x} < 0$) or ($\exp(-\Delta_{x'x}/T) > \text{random}[0, 1]$) THEN accept x′
    END DO

     IF *HEAT* THEN exit, ELSE calculate next annealing temperature T  
Until outer loops stopping criterion is satisfied  
where L is the length of Markov chain at T.

(2). Catch *messages on "Best So Far"*, inserting into inner and outer loops a procedure:

If improve "Best So Far", then index:=0;  
Otherwise, index:=index+1.

More in details, the procedure is:

*Memory-Annealing procedure*.  
x_best = x & f_best = f  
Repeat  
    DO k := 1 to L  
        Generate new solution x′  
        Calculate $\Delta_{x'x}$  
        IF ($\Delta_{x'x} < 0$) or ($\exp(-\Delta_{x'x}/T) > $ random[0, 1]) THEN  
          Accept x′  
          IF f < f_best THEN x_best = x & f_best = f  
        END IF  
    END DO  
    Calculate next annealing temperature T  
Until outer loops stopping criterion is satisfied  
x = x_best & f = f_best

(3). At the end of SA, take the optimal solution obtained as the initial solution to carry on a local search method or the SA method again.

*Annealing-Local Search procedure*.  
x_best = x & f_best = f; *anneal*=true  
Repeat  
    Repeat  
        DO k := 1 to L  
          Generate new solution x′ (Randomly local search or Allover local search)  
          Calculate $\Delta_{x'x}$  
          IF ($\Delta_{x'x} < 0$) or ( anneal & $\exp(-\Delta_{x'x}/T) >$ random[0, 1])) THEN  
            Accept x′  
            IF f < f_best THEN x_best = x & f_best = f  
          END IF  
        END DO  
        IF anneal THEN calculate next annealing temperature T  
    Until outer loops stopping criterion is satisfied  
    IF f_best < f THEN x = x_best & f = f_best  
    anneal:=not (anneal)  
Until anneal

In [52] after *Accept x′* a downhill local search method is embedded. On the contrary, before SA search we also may carry on a local search:

Jiapu Zhang

*Local Search-Annealing procedure.*
x_best = x & f_best = f; search=true & m := 0
Repeat
    IF search THEN L := Ls & m := m + 1 ELSE L := Lh
    Repeat
        a := 0
        DO k := 1 to L
            Generate new solution x′ (Random local search or All over local search)
            Calculate $\Delta_{x'x}$
            IF (search & $\Delta_{x'x} < 0$) or
            ( NOT search & $\Delta_{x'x} > 0$ & $\exp(-\Delta_{x'x}/T) >$ random[0, 1]) ) THEN
                Accept x′
                a := 1
            END IF
        END DO
    Until (search & (a = 0)) or (NOT search & (a = 1))
    IF f_best < f THEN x = x_best & f = f_best
    anneal:=not (anneal)
    search:=NOT search
Until m = snum
x = x_best & f = f_best
where snum is the number of optimal searches given.

(4). During the SA, for current state, take several search strategies, and accept the best state found with respect to probability.

## 3. The hybrid of SA with a local/global search method

Global optimization SA search sometimes trapped at local minima and cannot reach the real global minima. Hybrid with local search or global search optimization method is a strategy to bring the SA out of the trapped local minima. In this Section we introduce several successfully tested hybrid methods of SA. In Subsection 3.1, we will introduce the hybrid SA with local search discrete gradient (DG) method, and in Subsection 3.2 we will introduce the hybrids of SA with global search Self-Adaptive Evolutionary Strategy μ+λ (SAES(μ+λ)) method and global search Self-Adaptive Classical Evolutionary Programming (SACEP) method.

### 3.1 Hybrid SA discrete gradient method

#### 3.1.1 Efficiency of discrete gradient method

The DG method [53] is a derivative-free local search optimization method. Therefore, first we investigate the efficiency of discrete gradient method, comparing with other well-known derivative-free methods. We use small-size standard test problems to test the DG method, Nelder-Mead's simplex method [54], and Powell's UOBYQA method [54]. For each problem and each dimension, we run the three methods 50 times. The 50 initial solutions are randomly taken from the feasible region. The best optimal value obtained and its frequency of occurrence, the mean and the variance of 50 optimal values obtained can be seen in the following database.

*Numerical results for the DG method, Simplex method and UOBYQA method*

| Problem | Dimension | Method | Best value obtained | Frequency | Mean | Variance |
|---|---|---|---|---|---|---|
| Camel [55] | 2 | DG | -1.031628 | 80% | -0.86840 | 0.10876 |
| | | Simplex | -1.031628 | 82% | -0.85462 | 0.26388 |
| | | UOBYQA | -1.031628 | 46% | -0.17335 | 1.29922 |
| Goldstein-Price [56] | 2 | DG | 3.000000 | 50% | 65.10002 | 2.6E+04 |
| | | Simplex | 3.000000 | 60% | 372.30849 | 2.9E+06 |
| | | UOBYQA | 3.000000 | 40% | 122.88000 | 7.2E+04 |
| Griewanks [41] | 6 | DG | 9.224853 | 2% | 26.62373 | 121.67116 |
| | | Simplex | 0.277595 | 2% | 7.6E+04 | 2.9E+11 |
| | | UOBYQA | 0.946213 | 2% | 81.79783 | 1764.17313 |
| | 30 | DG | 3.233965 | 2% | 7.98272 | 5.45457 |
| | | Simplex | 67.904745 | 4% | 2.2E+06 | 3.3E+13 |
| | | UOBYQA | Failed | Failed | Failed | Failed |
| Hansen [57] | 2 | DG | -176.541793 | 94% | -174.67794 | 55.53534 |
| | | Simplex | -176.541793 | 44% | -134.99624 | 2652.23300 |
| | | UOBYQA | -176.541793 | 30% | -84.06059 | 5067.56707 |
| Hartman [58] | 3 | DG | -3.862782 | 100% | -3.86278 | 0.00000 |
| | | Simplex | -3.862782 | 76% | -3.71331 | 0.33637 |
| | | UOBYQA | -3.862782 | 80% | -3.49928 | 0.76195 |
| | 6 | DG | -3.322368 | 96% | -3.31760 | 0.00056 |
| | | Simplex | -3.322337 | 8% | -2.81436 | 0.32494 |
| | | UOBYQA | -3.322368 | 84% | -3.30329 | 0.00195 |
| Levy Nr.1 [59] | 2 | DG | 0.000000 | 100% | 0.00000 | 0.00000 |
| | | Simplex | 0.000000 | 52% | 1.80857 | 9.23844 |
| | | UOBYQA | 0.000000 | 34% | 4.19420 | 24.77323 |
| | 10 | DG | 0.000000 | 96% | 0.03110 | 0.03256 |
| | | Simplex | 0.025326 | 2% | 1.55279 | 2.71983 |
| | | UOBYQA | 0.000000 | 58% | 1.50279 | 8.14028 |
| | 20 | DG | 0.000000 | 96% | 0.00933 | 0.00238 |
| | | Simplex | 0.754527 | 2% | 4.63484 | 5.40594 |
| | | UOBYQA | 0.000000 | 70% | 1.07639 | 6.12405 |
| Levy Nr.2 [40] | 2 | DG | 0.000000 | 100% | 0.00000 | 0.00000 |
| | | Simplex | 0.000000 | 44% | 1.6E+07 | 1.9E+15 |
| | | UOBYQA | 0.000000 | 14% | 8.24957 | 47.28309 |
| | 10 | DG | 0.000000 | 86% | 0.04354 | 0.01188 |
| | | Simplex | 4.277692 | 2% | 2.1E+08 | 3.0E+16 |
| | | UOBYQA | 0.000000 | 6% | 1.69328 | 1.58222 |
| | 20 | DG | 0.000000 | 12% | 0.36077 | 0.06656 |
| | | Simplex | Failed | Failed | Failed | Failed |
| | | UOBYQA | 0.000000 | 2% | 1.24647 | 0.48673 |
| Levy Nr.3 [57] | 4 | DG | -21.502355 | 28% | -21.27131 | 0.16344 |
| | | Simplex | -21.499463 | 4% | -10.64044 | 108.39344 |
| | | UOBYQA | -18.392864 | 2% | 41.83276 | 1751.95111 |
| | 5 | DG | -11.504402 | 14% | -10.98174 | 0.25725 |
| | | Simplex | -11.489721 | 2% | -7.00003 | 23.97235 |
| | | UOBYQA | -8.505954 | 2% | 14.02895 | 210.25049 |
| Shekel-5 [58] | 4 | DG | -10.153199 | 96% | -9.85005 | 1.46911 |
| | | Simplex | -10.153200 | 16% | -6.11971 | 12.29406 |
| | | UOBYQA | -10.153200 | 74% | -8.44760 | 8.82110 |
| Shekel-7 [58] | 4 | DG | -10.402940 | 96% | -10.14444 | 1.69414 |
| | | Simplex | -10.402937 | 10% | -4.11194 | 8.81334 |
| | | UOBYQA | -10.402940 | 26% | -5.29144 | 9.83501 |
| Shekel-10 [58] | 4 | DG | -10.536410 | 96% | -10.26701 | 1.85441 |
| | | Simplex | -10.536313 | 4% | -3.39869 | 6.11953 |
| | | UOBYQA | -10.536410 | 20% | -4.45823 | 9.65360 |
| Shubert Nr.1 [58] | 2 | DG | -186.730908 | 92% | -179.82264 | 1.82815 |
| | | Simplex | -186.730909 | 52% | -132.55984 | 3915.58961 |
| | | UOBYQA | -186.730909 | 26% | -78.58946 | 5030.09156 |
| Shubert Nr.2 [58] | 2 | DG | -186.730908 | 18% | -168.91354 | 235.63292 |
| | | Simplex | -186.730909 | 2% | -115.19069 | 3686.41945 |
| | | UOBYQA | -186.730909 | 4% | -28.78545 | 5083.27183 |
| Shubert Nr.3 [57] | 2 | DG | -24.062499 | 42% | -22.59130 | 1.59932 |
| | | Simplex | -24.062499 | 32% | -19.43818 | 19.35178 |
| | | UOBYQA | -24.062499 | 12% | -14.99852 | 38.10882 |

By the comparative analysis from the above database, we find that the discrete gradient method is the best one among all those methods not only for low dimension problems but also for higher dimension problems. Nelder-Mead's simplex method cannot work well for the problems with dimensions greater than 10 and Powell's UOBYQA method cannot work fast for the problems with dimensions greater than 20. So, we choose the DG method.

The performance profile is a cumulative distribution function over a performance ratio and provides condensed information in terms of robustness, efficiency and quality of solution information. We briefly write the two formulae as follows:

$$\rho_{p,s}=t_{p,s}/\min\{t_{p,s}:1\leq s\leq n_s\} \text{ if } |(o_{p,s}-b_p)/b_p|\leq \delta, \text{ otherwise } \rho_{p,s}=\rho M,$$
$$p_s(\tau) =1/n_p \, C \, [\{p\in P: \rho_{p,s} \leq \tau\}]$$

where P is a given set of problems p=1, . . . , $n_p$, s is one of solvers s=1, . . . , $n_s$, $t_{p,s}$ is the solver resource (e.g. computational time) spent by solver s on problem p, $o_{p,s}$ denotes the solution value found by solver s for problem p, $b_p$ is the best solution value found when applying all solvers to problem p, $\delta>0$ is a user-defined relative objective value difference threshold, $\rho M$ is an upper bound on $\rho_{p,s}$ over all problem & solver pairs p, s in which solver s fails to solve problem p, $1\leq\tau\leq\infty$, and C [{·}] denotes the cardinality (size) of the set {·}. The function $\rho_{p,s}$ is called a performance ratio and the function $p_s(\tau)$ is called the performance profile function of the performance ratio. Our numerical experiments for the formula of $\rho_{p,s}$ are: (i) when $b_p$ is zero, in the denominator we replace bp by 1, and (ii) $\delta\leq10^{-4}$. The performance profile of the DG method, Simplex method and UOBYQA method showed that, being compared with the Simplex method and UOBYQA method, the DG method is absolutely the winner for solver resource and always better than the UOBYQA method for all the solver resource. This gave an explanation why we had chosen the DG method to use in this Subsection.

### 3.1.2 Hybrid DG-SA-DG algorithm

In this subsection we develop a hybrid SA and DG for solving the global optimization problem
$$\text{Minimize } f(x) \quad \text{subject to } x\in X,$$
where X is a subset of $R^n$, X is a compact set and f is continuous and it is also a locally Lipschitz continuous function. The hybrid method starts from an initial point, first executes the DG method to find local minimum, then carries on with the SA method in order to escape from this local minimum and to find a new starting point for DG method. Then we again apply the discrete gradient method starting from the current best point and so on until the sequence of the optimal objective function values obtained is convergent. The pseudo-code of the hybrid method is listed as following:

*Algorithm: Hybrid discrete gradient and simulated annealing method.*
**Initialization**:
Define the objective function f and its feasible solution space.
Call initial feasible solution generating procedure to get x.
Call initial temperature selecting procedure to get $T_0$.
Initialize f: f = f(x).
Initialize the neighborhood feasible solution x neighbour = 0.

Initialization of x_best: x_best = x.
Initialization of f_best: f_best = f.
do {
    **DG local search part**:
    F_best_local = local_search(x_best, x_new);
    x = x_new;
    **SA global search part**:
    do {
        do {
            x_neighbour = randomly_perturb(x);
            f_neighbour = f(x_neighbour);
            Calculate the difference $\Delta$=f_neighbour − f;
            If ($\Delta \leq 0$) or (random[0,1] < exp(-$\Delta$/Temperature))
                x = x_neighbour f = f_neighbour;
                If (f $\leq$ f_best) x_best = x f_best = f;
        } while (equilibrium has not been reached);
        Temperature annealing
    } while (Simulated Annealing stop criterion has not been met);
} while (f_best − f_best_local $\leq$ −0.001);

The convergence of the proposed hybrid method directly follows from the convergence of SA method and DG method. Generally, local search DG method makes the objective function value decrease a little bit from the initial guess, then the global search SA makes the value a big decreasing, the iterations go on until both the local and global searches cannot change the objective function value very much.

### 3.1.3 Implementations

**The Description of problems**
To examine the performance of Algorithm 14, we apply it to solve the wellknown complicated problems of Ackley, Bohachevsky, Branin, De Joung, Easom, Goldstein and Price, Griewank, Hartman, Hump, Hyper-Ellipsoid, Levy Nr.1, Levy Nr.2, Michalewicz, Neumaier Nr.2, Neumaier Nr.3, Rastringins, Rosenbrock, Schaffer Nr.1, Schaffer Nr.2, Shekel-N, Shubert Nr.1, Shubert Nr.2, Sphare, Step, Zakharov, Zimmermanns, which can be found from [40-41, 56, 58-59].

**Implementation of algorithm**
For the SA part of this hybrid method, we use the Neighborhood solution search procedure described at the end of Subsection 2.4. We still use T = 0.9*T as the cooling schedule in this hybrid method. The initial temperature is taken large enough according to the rule in [3]. The number of inner and outer iterations are taken to be large enough guaranteeing sufficient iterations. The DG method used here reduces the constrained minimization problem to unconstrained using exact penalty functions, and it terminates when the distance between the approximation to the subdifferential and origin is less than a given tolerance $\varepsilon > 0$ ($\varepsilon=10^{-4}$). The initial solution for the hybrid method is randomly taken from the feasible region of the problem.

**Results of numerical experiments and discussions**

Numerical experiments have been carried out in VPAC (Victorian Partnership for Advanced Computing) with CPU 833MHz. The results of numerical experiments are listed as follows. We denote HDGSAM the hybrid DG and SA method.

*Numerical results for HDGSAM*

| Function | Dimension | Best value obtained | Best value known | Number of function evaluations |
|---|---|---|---|---|
| Ackleys | 2 | 0.00008 | 0 | 200064 |
|  | 3 | 0.00012 | 0 | 200048 |
|  | 5 | 0.00009 | 0 | 200056 |
|  | 7 | 0.00008 | 0 | 200039 |
|  | 10 | 0.00003 | 0 | 200043 |
|  | 20 | 0.00052 | 0 | 2000623 |
|  | 30 | 0.00018 | 0 | 7000725 |
| Bohachevsky 1 | 2 | 0.000000 | 0 | 500020 |
| Bohachevsky 2 | 2 | 0.000000 | 0 | 500018 |
| Bohachevsky 3 | 2 | 0.000000 | 0 | 500017 |
| Branin | 2 | 0.397887 | 0.397887 | 300018 |
| De Joung | 3 | 0.0000000 | 0 | 500017 |
| Easom | 2 | -1.000000 | -1 | 20980 |
| Goldstein Price | 2 | 3.0000000 | 3 | 200016 |
| Griewank | 1 | 0.0000000 | 0 | 2068 |
|  | 2 | 0.0000000 | 0 | 2071 |
|  | 3 | 0.0000000 | 0 | 2071 |
|  | 4 | 0.0000000 | 0 | 2071 |
|  | 5 | 0.0000001 | 0 | 2072 |
|  | 6 | 0.0000001 | 0 | 6162 |
|  | 10 | 0.0000002 | 0 | 2085 |
|  | 20 | 0.0000003 | 0 | 2080 |
|  | 30 | 0.0000038 | 0 | 2088 |
| Hartmann | 3 | -3.86278215 | -3.86278215 | 6000091 |
|  | 6 | -3.3223680 | -3.3223680 | 6000165 |
| Hump | 2 | 0.000000 | 0 | 200030 |
| Hyper-Ellipsoid | 30 | 0.00000 | 0 | 3207 |
| Levy 1 | 5 | 0.000000 | 0 | 10006 |
|  | 10 | 0.000000 | 0 | 20379 |
|  | 20 | 0.000000 | 0 | 40971 |
|  | 30 | 0.000000 | 0 | 60661 |
|  | 50 | 0.000000 | 0 | 84454 |
|  | 70 | 0.0000000 | 0 | 140191 |
|  | 80 | 0.000000 | 0 | 140256 |
|  | 90 | 0.000000 | 0 | 180861 |
|  | 100 | 0.000000 | 0 | 203011 |
|  | 200 | 0.000000 | 0 | 400022 |
|  | 300 | 0.000000 | 0 | 604256 |
|  | 400 | 0.000000 | 0 | 805165 |
|  | 500 | 0.0000000 | 0 | 1002658 |
|  | 1000 | 0.000000 | 0 | 2025615 |
|  | 2000 | 0.000001 | 0 | 4064831 |
| Levy 2 | 5 | 0.000000 | 0 | 10089 |
|  | 10 | 0.000000 | 0 | 20225 |
|  | 20 | 0.000000 | 0 | 34429 |
|  | 30 | 0.0000000 | 0 | 60910 |
|  | 50 | 0.0000000 | 0 | 103288 |
|  | 70 | 0.0000000 | 0 | 96231 |
|  | 80 | 0.000000 | 0 | 162829 |
|  | 90 | 0.000000 | 0 | 123735 |
|  | 100 | 0.000000 | 0 | 71485 |
|  | 200 | 0.000000 | 0 | 148402 |
|  | 300 | 0.000000 | 0 | 603626 |
|  | 400 | 0.000000 | 0 | 283107 |
|  | 500 | 0.0000000 | 0 | 866036 |
|  | 1000 | 0.000000 | 0 | 533608 |
|  | 2000 | 0.000000 | 0 | 1948076 |
| Michalewicz | 2 | -1.8013 | -1.8013 | 20066 |
| Neumaier 2 | 4 | 0 | 0 | 6000369 |
| Neumaier 3 | 10 | -210 | -210 | 6094 |
| Rastringins | 2 | 0.0000000 | 0 | 200017 |
| | 3 | 0.0000000 | 0 | 200023 |
| | 5 | 0.0000000 | 0 | 200016 |
| | 7 | 0.0000000 | 0 | 200026 |
| | 10 | 0.0000000 | 0 | 200014 |
| Rosenbrock | 2 | 0.000000 | 0 | 1092 |
| | 5 | 0.000000 | 0 | 2086 |
| | 10 | 0.000000 | 0 | 2075 |
| Schaffer 1 | 2 | 0.0000026 | 0 | 2009038 |
| Schaffer 2 | 2 | 0.0000000 | 0 | 200015 |
| Shekel (N=5) | 4 | -10.15320 | -10.15320 | 2017770 |
| (N=7) | 4 | -10.40294 | -10.40294 | 6000159 |
| (N=10) | 4 | -10.53641 | -10.53641 | 200138 |
| Shubert 1 | 2 | -186.7309088 | -186.7309088 | 20000043 |
| Shubert 2 | 2 | -186.730903 | -186.730909 | 100025 |
| Sphare | 3 | 0.000000 | 0 | 6026 |
| Step | 5 | 0.000000 | 0 | 2134 |
| | 10 | 0.000000 | 0 | 2134 |
| | 50 | 0.000000 | 0 | 4204 |
| Zakharov | 2 | 0.00000000 | 0 | 300067 |
| | 5 | 0.00000000 | 0 | 300070 |
| | 10 | 0.00000000 | 0 | 300060 |
| Zimmermanns | 2 | 0.0000365 | 0 | 10000078 |

We can see the best objective function value obtained and the best objective function value known are equal to each other for every problem. This means our hybrid method is good and accurate for all those well-known optimization problems. Regarding the computational CPU time of the hybrid method for solving all those problems, it is very satisfactory. Take the Levy Nr.2 function as an example, the optimization method with 100 variables and 100 constraints needs 33.42019200 seconds to reach its optimal value 0.00000001, with 1000 variables and 1000 constraints needs 584.27654400 seconds to reach its optimal value 0.00000041, and even for the optimization problem with 10000 variables and 10000 constraints it needs 13696.099664 seconds to reach an objective function value 0.52427744. The runtime is on exponential increase. Based on results of numerical experiments, we can conclude for the hybrid DG and SA method: the hybrid DG and SA method is effective for many well-known optimization problems.

### 3.2 Hybrid SA evolutionary algorithms

In this Subsection, we present two hybrid SA and evolutionary algorithms.

Numerical results show that all these hybrid methods of evolutionary computation algorithms work well. The numerical results show us SAES($\mu+\lambda$) method [60] and SACEP method [60] can successfully work for all our test problems. The SA algorithm is a sequential computing algorithm and evolutionary algorithms are parallel computing algorithms. So, in this Subsection, using SA method, we improve them. We use SA as a search operator once for SAES($\mu + \lambda$) method, and once for SACEP method. Both the algorithms designed in this section simply work by applying the SA on all individuals in the population of the initial generation. In subsequent generations, SA is applied only for the best solutions found so far.

*Algorithm. SA-SAES($\mu + \lambda$).*

Step 0. Randomly generate μ parents, where each parent $z_k=(x_k, \sigma_k)$.
Step 1. Apply SA on each parent $x_k$.
Step 2. Set $\tau=\text{sqrt}(2\text{sqrt}(n))^{-1}$ and $\tau'=(\text{sqrt}(2n))^{-1}$.
Step 3. Until λ children are generated, do
Step 4. Select two parents $z_k=(x_k, \sigma_k)$ and $z_l=(x_l, \sigma_l)$ at random to generate child $y_j=(x_j, \sigma_j)$.
Step 5. Discrete recombination: for each variable $x_{ji}$ and step size $\sigma_{ji}$ in $y_j$, do ($x_{ji}=x_{ki}$ and $\sigma_{ji}=\sigma_{ki}$) or ($x_{ji}=x_{li}$ and $\sigma_{ji}=\sigma_{li}$)
Step 6. Mutation: For each $x_{ji}$ and step size $\sigma_{ji}$ in $y_j$
$$x'_{ji}=x_{ji}+\sigma_{ji}N_j(0, 1)$$
$$\sigma'_{ji}=\sigma_{ji}\exp(\tau'N(0, 1)+\tau N_j(0, 1))$$
Step 7. If the number of children is less than λ, go to Step 4.
Step 8. Select the best μ individuals among all the μ + λ parents and children.
Step 9. Apply SA on the best individual among the selected μ individuals.
Step 10. If the stopping criteria are satisfied, stop, else go to step 2.

*Algorithm. SA-SACEP.*
Step 0. Randomly generate μ parents and evaluate them, where each parent $z_k=(x_k, \sigma_k)$.
Step 1. Apply SA on each parent $x_k$.
Step 2. Set $\tau=\text{sqrt}(2\text{sqrt}(n))^{-1}$ and $\tau'=(\text{sqrt}(2n))^{-1}$.
Step 3. For each parent, generate a child as follows
$$x'_{ji}=x_{ji}+\sigma_{ji}N_j(0, 1)$$
$$\sigma'_{ji}=\sigma_{ji}\exp(\tau'N(0, 1)+\tau N_j(0, 1))$$
Step 4. Evaluate all children
Step 5. Undertake a tournament y for each parent and child as follows: select ζ individuals with replacement from the joint set of parents and children. For each individual z of the ζ individuals, if y is better than z, add 1 to the fitness of y.
Step 6. Select the best μ individuals among all parents and children with the highest fitness.
Step 7. Apply SA on the best individual among the selected μ individuals.
Step 8. If the stopping criteria are satisfied, stop, else go to step 1.

Numerical results listed as follows show that, from a point of view of the optimal objective function values obtained, the algorithms presented in this section separately improve SAES(μ + λ) method and SACEP method greatly.

*The Optimal objective function values of SAES(μ + λ) Algorithm and SA-SAES(μ + λ) Algorithm, and SACEP Algorithm and SA-SACEP Algorithm*

| Function | Number of variables | SAES(μ+λ) | SA-SAES(μ+λ) | SACEP | SA-SACEP |
|---|---|---|---|---|---|
| F1 [61] | 2 | -186.731 | -186.731 | -186.731 | -186.731 |
| F2 [62] | 5 | 1.0 | 1.0 | 1.0 | 1.0 |
|  | 10 | 1.0 | 1.0 | 1.0 | 1.0 |
|  | 20 | 1.28551 | 1.0 | 24.5297 | 1.0 |
|  | 30 | 1.02754 | 1.0 | 1.13336 | 1.0 |
|  | 50 | 1.00388 | 1.00001 | 9.28671 | 1.00001 |
| F3 (Ackleys) | 2 | 0.0 | 0.0 | 0.0 | 0.0 |
|  | 3 | 0.0 | 0.0 | 0.0 | 0.0 |
|  | 5 | 0.0 | 2.41563e-05 | 0.0 | 2.41563e-05 |
|  | 7 | 2.13384 | 4.86888e-05 | 1.72382 | 4.86888e-05 |
|  | 10 | 3.90647 | 7.6222e-05 | 1.08046 | 8.82517e-05 |
|  | 20 | 5.1886 | 0.000190629 | 2.24666 | 0.000224306 |
|  | 30 | 5.47366 | 0.0003507 | 4.92406 | 0.000406911 |
| F4 (Bohachevsky Nr.1) | 2 | 0.11754 | 0.117535 | 0.117548 | 0.117535 |
| F5 (Bohachevsky Nr.2) | 2 | 0.0 | 0.0 | 0.0 | 0.0 |
| F6 (Bohachevsky Nr.3) | 2 | 0.0 | 0.0 | 0.0 | 0.0 |
| F7 (Branin) | 2 | 0.398891 | 0.397887 | 0.398055 | 0.397887 |
| F8 (De Joung) | 3 | 0.0 | 0.0 | 0.0 | 0.0 |

| Function | n | | | | |
|---|---|---|---|---|---|
| F9 (Easom) | 2 | -0.999725 | -1.0 | -0.98863 | -1.0 |
| F10 (Goldstein Price) | 2 | 3.00006 | 3.0 | 3.0002 | 3.0 |
| F11 (Hartman with n = 3) | 3 | -3.86271 | -3.86278 | -3.86277 | -3.86278 |
| F12 (Hartman with n = 6) | 6 | -1.84847 | -3.32237 | -3.32192 | -3.32237 |
| F13 (Hump) | 2 | 8.86897e-05 | 4.65327e-08 | 0.000439177 | 0.0 |
| F14 (Hyper-Ellipsoid) | 30 | 1697.83 | 4.20078e-06 | 1.76103 | 0.0 |
| F15 (Levy Nr.2) | 5 | 0.0257144 | 1.02076e-10 | 0.0120519 | 0.0 |
| | 10 | 0.0129742 | 9.06744e-10 | 0.0317808 | 0.0 |
| | 20 | 2.34247e-06 | 5.48692e-09 | 0.0136671 | 0.0 |
| | 30 | 0.00193177 | 2.12137e-08 | 0.785024 | 0.0 |
| | 50 | 0.616365 | 6.12211e-08 | 2.07428 | 0.0 |
| F16 (Levy Nr.3) | 5 | 0.0218405 | 3.8796e-09 | 0.000743298 | 0.0 |
| | 10 | 0.00617594 | 1.35077e-08 | 0.000173664 | 0.0 |
| | 20 | 0.0 | 1.28154e-07 | 0.00358961 | 0.0 |
| | 30 | 0.000140932 | 4.54418e-07 | 0.000992482 | 0.0 |
| | 50 | 1.20497 | 1.68793e-06 | 1.32839e+06 | 1.60169e-06 |
| F17 (Michalewicsz) | 2 | -1.95063 | -1.8013 | -1.95217 | -1.8013 |
| F18 (Neumaier Nr.2) | 4 | 0.00245258 | 0.000487242 | 0.0766711 | 0.000174267 |
| F19 (Neumaier Nr.3) | 10 | -21.0244 | -209.998 | -203.925 | -209.999 |
| F20 (Rastringins Nr.1) | 2 | 0.0 | 2.36476e-10 | 0.0 | 0.0 |
| | 3 | 0.0 | 3.91857e-10 | 0.995047 | 0.0 |
| | 5 | 0.0 | 3.25394e-08 | 5.97189 | 0.0 |
| | 7 | 0.0 | 1.93565e-07 | 8.95636 | 0.0 |
| | 10 | 1.99124 | 1.98263e-06 | 32.8386 | 1.98263e-06 |
| F21 (Rosenbrock) | 2 | 0.0079492 | 5.68257e-07 | 0.00856004 | 2.096e-06 |
| | 5 | 0.915901 | 0.000190216 | 0.00588099 | 3.13482e-05 |
| | 10 | 4.104 | 3.83856e-05 | 2.15272 | 0.000239605 |
| F22 (Schaffer Nr.1) | 2 | 0.0 | 0.0 | 0.0 | 0.0 |
| F23 (Schaffer Nr.2) | 2 | 0.0 | 0.195296 | 0.0 | 0.195296 |
| F24 (Shekel-5) | 4 | -5.04985 | -5.27766e+13 | -5.05082 | -5.27766e+13 |
| F25 (Shekel-7) | 4 | -5.0606 | -5.27766e+13 | -5.05484 | -5.27766e+13 |
| F26 (Shekel-10) | 4 | -5.1273 | -5.27766e+13 | -5.11435 | -5.27766e+13 |
| F27 (Shubert Nr.1) | 2 | -186.731 | -186.731 | -186.731 | -186.731 |
| F28 (Shubert Nr.2) | 2 | -186.341 | -186.731 | -186.731 | -186.731 |
| F29 (Step) | 5 | -144.0 | 0.0 | -2848 | 0.0 |
| | 10 | -366 | 0.0 | -1.18937e+07 | 0.0 |
| | 50 | -13864 | 0.0 | -7.19852e+34 | 0.0 |
| F31 (Zimmermanns) | 2 | -103.806 | -494.741 | -494.748 | -494.735 |
| F32 [58] | 2 | 0.0 | 4.19095e-06 | 0.0 | 4.19095e-06 |
| | 5 | 0.0 | 6.11739e-05 | 0.0 | 6.11739e-05 |
| | 10 | 0.0 | 0.000461433 | 0.00287121 | 0.000553783 |
| | 50 | 0.681216 | 6.26669e-13 | 15.1833 | 0.017286 |
| F33 [58] | 2 | 0.0 | 6.26669e-13 | 0.0 | 0.0 |
| | 5 | 0.0 | 4.16862e-06 | 0.0 | 4.16862e-06 |
| | 10 | 8.06556 | 0.0113471 | 0.0 | 0.0322543 |
| | 50 | 11206.3 | 1030.77 | 6733.64 | 894.608 |
| F34 [58] | 2 | 0.0 | 1.78348e-05 | 0.0 | 1.78348e-05 |
| | 5 | 0.0 | 0.000742218 | 0.0 | 0.000742218 |
| | 10 | 0.0 | 0.00371919 | 0.0 | 0.00371919 |
| | 50 | 79.0741 | 0.189146 | 14.1295 | 0.189135 |
| F35 [58] | 2 | 0.0 | 0.0 | 0.0 | 0.0 |
| | 5 | 0.0 | 0.0 | 0.0 | 0.0 |
| | 10 | 0.0 | 0.0 | 0.0 | 0.0 |
| | 50 | 52.0 | 0.0 | 33.0 | 0.0 |
| F36 [58] | 2 | 5.03179e-06 | 8.18303e-07 | 4.59426e-06 | 8.18303e-07 |
| | 5 | 2.03186e-05 | 1.28561e-05 | 0.000226618 | 6.52515e-06 |
| | 10 | 0.000277681 | 6.09379e-05 | 0.001168 | 6.73284e-05 |
| | 50 | 415.836 | 0.00306063 | 380.029 | 0.00440367 |
| F37 [58] | 2 | -837.931 | -837.966 | -3947.21 | -837.966 |
| | 5 | -1796.66 | -2094.91 | -1513.87 | -2094.91 |
| | 10 | -3809.75 | -4189.83 | -3245.56 | -4189.83 |
| | 50 | -18813.3 | -20949.1 | -12809.6 | -20949.1 |
| F38 [58] | 2 | 0.0179898 | 6.08096e-11 | 0.00202397 | 0.0 |
| | 5 | 0.0744221 | 9.45082e-09 | 0.0409532 | 0.0 |
| | 10 | 0.0019571 | 1.13757e-07 | 0.0114677 | 0.0 |
| | 50 | 7.73384 | 4.46473e-05 | 10.5956 | 4.02717e-05 |
| F41 [58] | 2 | -4.12397 | -4.12398 | -4.12373 | -4.12398 |
| F42 [58] | 2 | 0.398891 | 0.397887 | 0.398055 | 0.397887 |
| F43 [58] | 2 | 3.00006 | 3.0 | 3.00002 | 3.0 |

For problems F2~F4, SA-SAES($\mu+\lambda$) and SA-SACEP perform similarly. For problems F5~F18, all these two hybrid methods seem to have the same performance too. For problems F19~F22, SA-SAES($\mu+\lambda$) and SA-SACEP perform well. For F23~F26, SA-SAES($\mu+\lambda$) and SA-SACEP perform the same. For the problems remained, we cannot make some comparison for all the two hybrid methods of this Subsection. However, all the two hybrid methods in this Subsection are better than without hybrids.

## 4. Applying SA to molecular modeling

In [62], Zhang (2011) used the follow SA procedure: "The solvated proteins were then quickly heated from 0 K to 300 K linearly for 20 ps. The systems were kept at 300 K for 80 ps. The systems then were slowly cooled from 300 K to 100 K linearly for 400 ps. The systems were kept for 100 ps at 100 K. All the systems were in constant NVT ensembles using Langevin thermostat algorithm with weak restraints (a force constant of 10.0 kcal mol$^{-1}$ Å$^{-2}$ was used) on the solvated proteins. The SHAKE and SANDER (simulated annealing with NMR-derived energy restraints) algorithms with nonbonded cutoffs of 9 Å were used during the heating, cooling and the 100 ps at 100 K. The equilibration was done in constant NPT ensembles under a Langevin thermostat for 4,400 ps and the RMSD, PRESS, and VOLUME (DENSITY) were sufficiently stable for each model; the jump in RMSD of around 0.2 Å correlates with removing restraints for the change from NVT to NPT, but it did not change the structures at 100 K. Equilibration was under constant pressure 1 atm and constant temperature (100 K) in a neutral pH environment (equilibration was performed at the low temperature of 100 K in order to be consistent with the experimental work). A step size of 2 fs was used for all SA runs. The structures were saved to file every 100 steps. During the SA, the Metropolis criterion was used." This is an ordinary SA method for molecular modeling and refinement of 100K-crystal structures. Before or after the SA procedure is always using hybrids with local search optimization methods (for example the steepest descent method or the conjugate gradient method or both) in Amber and Gromacs computational chemistry packages.

## 5. Conclusions

SA is a popular method used in mathematical optimization computations. This paper introduced the detailed implementations of SA in optimization, and presented two kinds of hybrids with local and global optimization search algorithms respectively, and then introduced the situation in the use of SA in molecular modeling in computational chemistry crystal structures. The SA theory presented in this paper should be very useful for the practical SA calculations.

**Acknowledgements**

Dr Zhang JP thanks Professors Bagirov AM and Abbass HA for their computer codes of DG method and evolutionary algorithms.